\documentclass[12pt]{amsart}
\usepackage{amssymb}
\begin{document}
\numberwithin{equation}{section}

\def\1#1{\overline{#1}}
\def\2#1{\widetilde{#1}}
\def\3#1{\widehat{#1}}
\def\4#1{\mathbb{#1}}
\def\5#1{\frak{#1}}
\def\6#1{{\mathcal{#1}}}

\def\C{{\4C}}
\def\R{{\4R}}
\def\N{{\4N}}
\def\Z{{\4Z}}
\def\Q{{\4Q}}

\title[Deformation of generic submanifolds]{Deformation of
generic submanifolds in a complex manifold}
\author[M. S. Baouendi, L. P. Rothschild,  D. Zaitsev]{M. S.
Baouendi, L. P. Rothschild, D. Zaitsev}
\address{ Department of Mathematics, University of
California at San Diego, La Jolla, CA 92093-0112}
\email{sbaouendi@ucsd.edu,
lrothschild@ucsd.edu }
\address{ School of Mathematics, Trinity College, Dublin 2,
IRELAND}
\email{zaitsev@maths.tcd.ie }

\thanks{
The research of the first and second authors were supported in part by the NSF grant DMS-0400880.
The research of the third author was supported in part by the RCBS grant of the Trinity College Dublin
and by the Science Foundation Ireland grant 06/RFP/MAT018}
\thanks{ 2000 {\it   Mathematics Subject Classification.}  32H35, 32V40, 57A35, 58A20}
\thanks{{\it Keywords:} generic manifold, complex manifold, nondegenerate, finite type, deformation.}
\maketitle

\begin{abstract}
This paper shows that 
 an arbitrary generic submanifold in a complex
manifold can be deformed into a 1-parameter family of generic submanifolds satisfying strong nondegeneracy conditions. The proofs use a careful analysis of the jet spaces of embeddings satisfying certain nondegeneracy properties, and also make use of the Thom transversality   theorem, as well as the stratification of real-algebraic sets.  Optimal results on the order of nondegeneracy are given.
\end{abstract}

\def\Label#1{\label{#1}}

\def\cn{{\C^n}}
\def\cnn{{\C^{n'}}}
\def\ocn{\2{\C^n}}
\def\ocnn{\2{\C^{n'}}}


\def\dist{{\rm dist}}
\def\const{{\rm const}}
\def\rk{{\rm rank\,}}
\def\id{{\sf id}}
\def\aut{{\sf aut}}
\def\Aut{{\sf Aut}}
\def\Diff{{\sf Diff}}
\def\CR{{\rm CR}}
\def\GL{{\sf GL}}
\def\Re{{\sf Re}\,}
\def\Im{{\sf Im}\,}

\def\codim{{\rm codim}}
\def\crd{\dim_{{\rm CR}}}
\def\crc{{\rm codim_{CR}}}

\def\phi{\varphi}
\def\eps{\varepsilon}
\def\d{\partial}
\def\a{\alpha}
\def\b{\beta}
\def\g{\gamma}
\def\G{\Gamma}
\def\D{\Delta}
\def\Om{\Omega}
\def\k{\kappa}
\def\l{\lambda}
\def\L{\Lambda}
\def\z{{\bar z}}
\def\w{{\bar w}}
\def\t{\tau}
\def\th{\theta}

\newcommand{\rank}{\mathop{{\rm rank}}}
\newcommand{\tensor}{\otimes}
\newcommand{\Hom}{\mathop{Hom}}

\emergencystretch15pt
\frenchspacing

\newtheorem{Thm}{Theorem}[section]
\newtheorem{Cor}[Thm]{Corollary}
\newtheorem{Pro}[Thm]{Proposition}
\newtheorem{Lem}[Thm]{Lemma}

\theoremstyle{definition}\newtheorem{Def}[Thm]{Definition}

\theoremstyle{remark}
\newtheorem{Rem}[Thm]{Remark}
\newtheorem{Exa}[Thm]{Example}
\newtheorem{Exs}[Thm]{Examples}

\def\bl{\begin{Lem}}
\def\el{\end{Lem}}
\def\bp{\begin{Pro}}
\def\ep{\end{Pro}}
\def\bt{\begin{Thm}}
\def\et{\end{Thm}}
\def\bc{\begin{Cor}}
\def\ec{\end{Cor}}
\def\bd{\begin{Def}}
\def\ed{\end{Def}}
\def\br{\begin{Rem}}
\def\er{\end{Rem}}
\def\be{\begin{Exa}}
\def\ee{\end{Exa}}
\def\bpf{\begin{proof}}
\def\epf{\end{proof}}
\def\ben{\begin{enumerate}}
\def\een{\end{enumerate}}
\def\beq{\begin{equation}}
\def\eeq{\end{equation}}

\section{Introduction} Many interesting properties of real
submanifolds in complex space or, more generally, in a
complex manifold, require the imposition of some
nondegeneracy conditions.  For hypersurfaces, and in particular
for boundaries of domains in
$\C^N$,  the most studied condition is that of
Levi-nondegeneracy.  Although it is well-known that the
boundary of any weakly pseudoconvex domain can be
deformed into a strictly pseudoconvex one, it is not true that an
arbitrary smooth boundary  can be deformed into a (smooth)
boundary that is everywhere Levi-nondegenerate. In contrast,
in this paper we show that any (not totally real) generic
submanifold of
a complex manifold can be deformed into another generic
submanifold satisfying a higher order nondegeneracy condition
at every point of the deformed submanifold. The higher order
nondegeneracy conditions considered here have been used to
establish many properties that were previously only known
to hold  for Levi-nondegenerate hypersurfaces.

The first nondegeneracy condition we consider here is that of
finite nondegeneracy or, more precisely,
$k$-nondegeneracy for some integer $k\ge 1$, generalizing
Levi nondegeneracy. In particular, a
hypersurface $M\subset\C^N$ is
$1$-nondegenerate at a point $p$ if and only if it is
Levi-nondegenerate at $p$. The second nondegeneracy
condition we consider  is that of finite type in the sense of
Kohn \cite{Kohn} and Bloom-Graham \cite{BG}. For
hypersurfaces, finite nondegeneracy implies the finite type
property, but the two notions are independent for generic
submanifolds of higher codimension. The precise definitions of
these two conditions will be given in Section~\ref{def}.
The importance of these conditions is illustrated, for instance,
by the fact that, taken together, they are sufficient to guarantee
that the full CR automorphism group of a generic
submanifold, equipped with its natural topology,  is  a
finite-dimensional Lie group
   (see \cite{BRWZ}). Also, finite
nondegeneracy and finite type at a given point $p$ of a
real-analytic generic  submanifold $M\subset\C^N$ imply
that the stability group of $M$ at $p$ (i.e.\ the group of germs
of biholomorphisms fixing $p$ and preserving $M$) is a Lie
group, see
\cite{CM74, BS, Bel, BERasian,Z97,BERrat}.

The present paper is devoted to
approximation  of general generic submanifolds of a complex
manifold (with respect to the
strong or Whitney topology, see e.g.\
\cite{GG,Hir}) by those satisfying the  two nondegeneracy
conditions mentioned above.  More precisely,  our first result is
the following.

\bt\Label{weak} Let $M$ be a smooth (resp. real-analytic)
connected manifold  and $X$ be a complex manifold
such that $\dim_\C X<\dim_\R M<2 \dim_\C X$. Let
$\tau\colon M\to X$ be  a
smooth (resp. real-analytic) embedding such that
$\tau(M)$ is a generic submanifold of $X$. Then
any open neighborhood of $\tau$ in $C^\infty(M,X)$,  equipped with
the Whitney topology, contains a smooth (resp. real-analytic)
embedding whose image is a generic submanifold that is both
finitely nondegenerate and of finite type.
\et

For the case of a smooth, generic embedding, we shall also
prove a stronger version of  Theorem~\ref{weak} in which the
approximation is accomplished via  a
$1$-parameter family of embeddings. We have the
following.
\bt\Label{strong}
Let $M$ be a smooth
connected manifold, $X$ a complex manifold with $$\dim_\C
X<\dim_\R M<2\dim_\C X$$ and
$\tau\colon M\to X$ a smooth  embedding
such that $\tau(M)$ is a generic submanifold of
$X$. Then for any open neighborhood $U$ of $\tau$ in $C^\infty(M,X)$,
  there exists a smooth map
$\2\tau\colon M\times(-1,1)\to X$  such that
$\2\tau(p,0)=\tau(p)$ for all $p\in M$, and, for each $t\ne 0$,
the map $p\mapsto
\2\tau(p,t)$ is in $U$ and is an embedding of $M$ into $X$
whose image is a generic manifold that is finitely
nondegenerate and of finite type.
\et

In  the course of the proof of Theorems \ref{weak} and
\ref{strong}, we shall give explicit estimates on the order of the
nondegeneracy and the type of the approximating family of
embeddings.  This is made more
precise in Theorems \ref{weak-precise} and
\ref {strong-precise} below.  In fact, the nondegeneracy and the
type of the approximating embeddings can be estimated in
terms of the dimensions of $M$ and $X$ alone.  We will
also show that some of these estimates are sharp, as is
illustrated by Example~\ref{ge}; see also Remark~\ref{rem}.

A number of authors have studied deformations of
Levi-nondegenerate  real hypersurfaces and their embeddings into a
complex manifold.  Of the work in this direction, we mention here,
in particular, the papers of Burns-Shnider-Wells \cite{BSW},
Catlin-Lempert \cite{CL}, Bland-Epstein
\cite{BE}, Huisken-Klingenberg
\cite{HK}, and Huang-Luk-Yau \cite{HLY}.  

Two of our main tools in the proof of our results are the
Thom transversality theorem  and the stratification of
semi-algebraic sets (see e.g.\ \cite{BR}). The paper is organized as follows.  In Section
\ref{def} we give various definitions used throughout the paper and we
state Theorems \ref{weak-precise} and \ref{strong-precise}, which are
   the more precise versions of Theorems \ref{weak} and
\ref{strong}, giving in particular the order of nondegeneracy of the
approximating manifolds.  In Sections \ref{decomp} through
\ref{s:strong} we describe a precise stratification of jets of degenerate
embeddings and of those not of strong type (see Definition~\ref{str-type} below) and calculate their
codimension in the space of all jets.  The proof of Theorems
\ref{weak-precise} is given in Section \ref{proofweak} and that of
Theorem \ref{strong-precise} is given in Section \ref{proofstrong}. In
Section \ref{remarks}, we conclude with some remarks and examples.

\section{Definitions and precise results}\Label{def}
Let $\6M$ be a real submanifold in a complex manifold $X$.
Recall that $\6M$ is called {\em generic} if, for every $p\in
\6M$,
$$T_p\6M + JT_p\6M = T_pX,$$
where $J$ denotes the complex structure of $X$, and $T_p\6M$
and $T_pX$ are  the real tangent spaces of $\6M$
and $X$ respectively.  In what follows,
$\6M$ is assumed to be generic.  We let
$$N:=\dim_{\C}X, \ \ \ m:= \dim_{\R}\6M, \ \ \ d: =
\codim_{\R}\6M= 2N-m,$$
and denote
by $T^{(0,1)}\6M$ the bundle of $(0,1)$ vector fields on $\6M$.
  Then for any $p\in \6M$ we have
$\dim_{\C}T^{(0,1)}=N-d=:n$.  We note here that the
condition $N < m <2N$ imposed in the hypotheses of
Theorems \ref{weak} and \ref{strong} for $\6M = \tau(M)$ is
equivalent to the conditions $d >0$  and $n>0$.  We shall
assume these conditions in the remainder of  this paper.

We shall need the following definitions.
Given $p\in \6M$ and a system of local holomorphic coordinates
$Z=(Z_1,\ldots,Z_N)$ in $X$ near $p$, let
$\rho=(\rho^1,\ldots,\rho^d)$ be a system of smooth, real
defining functions of $\6M$ near $p$ (i.e.\
$\partial\rho^1(p,\1p)\wedge\cdots\wedge
\partial\rho^d(p,\1p)\ne0$).  For an integer $k\ge 1$, define the
{\em $k$-degeneracy} of
$\6M$ at $p$ as the integer
\begin{equation}\Label{rp}
r^1_p(k):=N-\dim{\sf span}_\C \left\{(L_1\ldots
L_s\rho^j_Z)(p,\1p) :  1\le j\le d;
0\le s \le k \right\},
\end{equation}
where $L_1,\ldots,L_s$ runs through arbitrary systems of
$(0,1)$ vector fields on $\6M$, i.e.\ sections of $T^{0,1}\6M$.
Here $\rho^j_Z=(\rho^j_{Z_1}, \ldots, \rho^j_{Z_N})$ denotes
the complex gradient of $\rho^j$ (with respect to the
coordinates $Z$) and is regarded as a vector in
$\C^N$. Observe that $0 \le r^1_p(k) \le N-d$.  It can be shown, see
e.g.
\cite[Chapter 11]{BER}, that the right-hand side in \eqref{rp} is an
invariant, i.e. independent of the choice of the defining function $\rho$
and the complex coordinates $Z$.

Recall that the generic submanifold $\6M\subset X$ is called
{\em $k$-nondegenerate at $p\in \6M$} if and only if $r^1_p(k)=0$
and $r^1_p(j)>0$ for $1\le j<k$.
Furthermore, $\6M$ is said to be {\em finitely nondegenerate at
$p$} if  it is $k$-nondegenerate at $p$ for some $k\ge 1$.
Finally, we say that $\6M$ is
{\em $k$-degenerate at $p$} if $r^1_p(k)>0$. With this
definition, a hypersurface $\6M\subset X$ is Levi
nondegenerate at $p$ if and only if it is 1-nondegenerate at $p$.

Recall also that the generic real submanifold $\6M\subset X$ is
said to be of {\em finite type $k$} at
$p\in \6M$ (in the sense of {\sc Bloom-Graham-Kohn}) if $\C
T_p\6M$ is spanned by the commutators of the form
\begin{equation}\Label{comm0}
[\6L_1,[\6L_2,\ldots,[\6L_{l-1},\6L_l]\ldots]](p),
\quad 1\le l\le k,
\end{equation}
where each $\6L_r$, $1\le r\le
l$, is either a $(1,0)$ or  a $(0,1)$ vector field on $\6M$,
and $k$ is minimal
with this property.
Here, for $l=1$, the quantity \eqref{comm0} is simply $\6L_1(p)$.
Furthermore, $\6M$ is said to be of {\em finite type at $p$} if it
is of finite type $k$ for some $k\ge1$.

We shall now introduce a stronger version of the above mentioned finite
type condition, which turns out to be more easily computable in terms of
the defining functions of the manifold $\6M$ (see Lemma
\ref{type} below).

\bd\Label{str-type}
Let $\6M\subset X$  be a  generic submanifold and $1\le
k<\infty$ be an integer.  We say that $\6M$ is of {\em strong type} $k$  at $p$  if
  $\C T_p\6M$ is spanned by the commutators of
the form
\begin{equation}\Label{comm}
[L_1,[L_2,\ldots,[L_{l-1},\1L_l]\ldots]](p),
\quad 1\le l\le k,
\end{equation}
and their complex conjugates, where $L_r$, $1\le r\le
l$, are any $(1,0)$ vector fields on $\6M$, and $k$ is minimal
with this property.
Here, for $l=1$, the quantity \eqref{comm} is simply $\1L_1(p)$.
\ed

It is an immediate consequence of the definition that, if
$M$ is of strong type $k$ at $p$, then it is of finite type $\le
k$ at $p$.
Moreover, for $k=2,3$, the notions of finite type $k$ and strong type $k$
coincide. However,  for
$k\ge 4$, finite type $k$ may not imply strong type $l$ for any $l$, even
for hypersurfaces. For example, the hypersurface $\Im w= |z|^4$ in
$\C^2$ is of finite type $4$ at $0$ but not of strong type $l$
for any $l\ge 1$.

Note that in the case of hypersurfaces in $X$,
$k$-nondegeneracy implies strong type $\le k+1$. On the
other  hand, in higher codimension, the conditions of being of
strong type $\le k$ and $l$-nondegeneracy are independent for any $k$
and $l$ (see Lemma
\ref{type} below).

For an arbitrary generic manifold $\6M\subset X$
and $p\in \6M$ and an integer $k \ge 1$, we define the {\em
$k$-defect of $\6M$ at $p$} as the nonnegative integer
$r^2_p(k)$ given by
\begin{equation}\Label{defect} r^2_p(k):=\dim_{\R} \6M-\dim_{\C}
V,\end{equation}
  where
$V$ is the span of the commutators in
\eqref{comm} and their complex conjugates.
Note that $\6M$ is of strong type $k$ at $p$ if and only if
$r^2_p(k)=0$ and $r^2_p(j)>0$ for $1\le j<k$.

Before stating a more precise version of Theorem
\ref{weak}, from which the latter follows, we introduce the following notation.
For any pair of positive integers $(m,N)$, with $2\le N < m< 2N$, we let
  $k_1(m,N)$ be the positive integer defined by
\begin{equation}\Label{k1} k_1(m,N) = \begin{cases}1\quad \text{if }\ N+2 \le
m\le 2N-3,  (m,N)\not= (7,5),
\cr 3 \quad \text{if } (m,N)=(3,2),\cr 2\quad \text{otherwise},\end{cases}
\end{equation}
and $k_2(m,N)$ be the smallest positive integer $k$ for which the following
inequality holds:
\begin{equation} \Label{k2} 2(m-N)\Big({k+m-N-1\atop k-1}\Big)
\ge (m-N)^2+2m.\end{equation}
Observe in particular that if $m = 2N-1$ (which will correspond to the case
of a hypersurface in Theorem \ref{weak-precise} below) then the integer
$k_2(m,N)$ is $4$ for $N=2$,
it is $3$ for $N=3$, and it is $2$ for $N>3$.
\bt\Label{weak-precise}
Let $M$, $X$ and $\tau\colon M\to X$ be as in Theorem~\ref{weak}
and set $m:=\dim_\R M$,
$N:=\dim_\C X$.  Let $k_1=k_1(m,N)$ and $k_2=k_2(m,N)$ be defined by
\eqref {k1} and \eqref{k2} respectively. Then  any neighborhood of $\tau$ in
$C^\infty(M,X)$, equipped with the Whitney topology,  contains a smooth (resp.
real-analytic)  embedding whose image is a generic submanifold that, at every
point, is both
$\ell_1$-nondegenerate, for some $\ell_1 \le
k_1$, and of strong type $\ell_2$, for some
$\ell_2 \le k_2$.
\et

\br
In contrast to $k_1$, the integer $k_2$ in
Theorem~\ref{weak-precise} has
no uniform bound for all $m$ and $N$.
Indeed, one can check by using Lemma~\ref{type} below that
   if $\6M\subset \C^N$ is a real $(N+1)$-dimensional
generic submanifold
of strong type $k$ at some point, then necessarily $N<2k$.
\er

In analogy with Theorem~\ref{weak-precise}, we have the following
more precise version of Theorem~\ref{strong}, from which
the latter follows.

\bt\Label{strong-precise}
Let $M$, $X$ and $\tau\colon M\to X$ be as in
Theorem~\ref{strong}. Then there exist  positive integers
$k'_1$ and $k'_2$,
depending only
on the dimensions
$m:=\dim_\R M$ and
$N:=\dim_\C X$,
such that
  any neighborhood $U$ of $\tau$ in $C^\infty(M,X)$
contains a smooth $1$-parameter
deformation $\2\tau$ of $\tau$ as in Theorem~\ref{strong}
such that for each $t\ne0$, the image of the embedding $p\mapsto
\2\tau(p,t)$ is, at each point $p$,    both
$\ell_1$-nondegenerate, for some $\ell_1 \le
k'_1$, and of strong type $\ell_2$, for some
$\ell_2 \le k'_2$.  Moreover, when $X=\C^N$,
$k'_1$ can be chosen to be the same as the integer $k_1(m,N)$ given by
\eqref{k1} and
$k'_2$ to  be the smallest positive integer $k$, for which the
following holds:
\begin{equation}\Label {k2'}
2(m-N)\Big({k+m-N-1\atop k-1}\Big)
\ge (m-N)^2+2m+1.
\end{equation}
In particular, if
$\tau(M)$ is a real hypersurface in $\C^N$,
i.e. $m=2N-1$, $k'_2$ can
  be chosen to be
$4$ for
$N=2$,
to be $3$ for
$N=3$, and
to be $2$ for $N>3$.
\et

The following is an immediate consequence of  Theorem
\ref{strong-precise}.
\bc Let $M\subset\C^N$ be a real hypersurface with $N\ge 3$.
Then $M$ can be smoothly approximated by a
$1$-parameter family of real hypersurfaces that are
$2$-nondegenerate  and of strong type $3$.
\ec

In  contrast, for $N=2$, there are real-analytic hypersurfaces
in $\C^2$ that cannot be approximated by
$2$-nondegenerate ones, as Example \ref{no2} below shows.

\section{Decomposition of jet spaces}\Label{decomp}

In what follows, if  $A$ and $B$ are smooth manifolds and $k$
a positive integer, we use the standard notation $J^k(A,B)$
for the manifold of all k-jets of smooth mappings from $A$ to
$B$. For
$p\in A$ and $q\in B$, we denote by $J_p^k(A,B)$ the
submanifold of all $k$-jets in $J^k(A,B)$ with source $p$ and by
$J_{p,q}^k(A,B)$ --- of those with source $p$ and target $q$. 
(See e.g.\ \cite{GG} for properties of these manifolds.)
If $f:A\to B$ is a map of class $C^k$ and $p\in A$, we denote by
$j^k_p f \in J^k_p(A,B)$ the jet of $f$ at $p$. 

Recall that if $f:A\to B$ is a smooth mapping, we say that $f$
is an embedding if $f$ is an immersion at every point of $A$
and is homeomorphism onto its image $f(A)$. If $M$ and $X$
are as in Theorem \ref{weak}, an embedding $f:M\to X$ is
called {\it generic} if $f(M)$ is a generic submanifold of $X$.
If $p\in M$ and $g:(M,p) \to X$ is a germ of a smooth mapping
at $p$, then $g$ is a {\it germ of a  generic embedding} if
$g'(p):T_pM \to T_{g(p)}X$ is injective and
$$T_{g(p)}X = g'(p) T_pM + J g'(p)T_pM.  $$
For $k \ge 1$ and  $p \in M$, we denote by $W_p^0\subset
J^k_p(M,X)$ the
$k$-jets of germs of generic embeddings as defined above.

If
$U$ is a local coordinate chart on $M$, then  it is standard
(see e.g.\ \cite{GG}) that
$J^k(U,\C^N)$ can be naturally identified with
$U\times \R^K$ for $K:=\dim_\R J^k_p(U,\C^N)$ for $p\in
U$.

Before stating the main result of this section, Lemma \ref{split}, we
need the following.

\bl\Label{germ} Let $Z=(Z_1,\ldots,Z_N)$ be fixed
coordinates in $\C^N$
and
$g:(\R^{2n+d},p)\to\C^{N}$  a germ of a generic embedding.
Then after reordering the coordinates $Z_j$ and multiplying some
by
$\sqrt {-1}$ if necessary, we can write $Z=(z,w)= (z_1,
\ldots,z_n,w_1,\ldots, w_d)$ such that the following holds.
After identifying $\C^N$ with $\R^n_x \times \R^n_y\times
\R^d_s\times \R^d_t$ by $(x,y,s,t)= (\Re z,\Im z, \Re w, \Im w)$
and writing $g=(g_1,g_2)$ with $g_1:(\R^{2n+d},p) \to \R^n_x
\times
\R^n_y\times
\R^d_s$ and $g_2:(\R^{2n+d},p) \to\R^d_t$, one has that $g_1$ is a
germ of a diffeomorphism at $p$, and for $U$ a sufficiently small
neighborhood of $p$ in $\R^{2n+d}$,
$g(U)$ is given near $g(p)$ by
\beq\Label{phi'}
\Im w=\phi(\Re z,\Im z,\Re w),
\end{equation} where $\phi$ is the germ   at
$g_1(p)$ given by
$\phi =g_2\circ g_1^{-1}$.  Moreover, the $d \times d$ matrix
$(\id + i {\partial \phi/ \partial s})$ is invertible at $g_1(p)$.
\el
\bpf  For  a sufficiently small
neighborhood $U$ of $p$ in $\R^{2n+d}$, let
$\rho=(\rho_1,\ldots,\rho_d)$ be some real-valued local defining
functions for $g(U)$ near $g(p)$.  Since the $N\times d$ matrix
$(\partial \rho/\partial Z)$ has rank $d$ at $g(p)$, after reordering
the coordinates $(Z_1,\ldots,Z_N)$ and multiplying some of the
$Z_j$ by $\sqrt {-1}$ if necessary, we may assume that
$Z=(z,w) = (z_1,
\ldots,z_n,w_1,\ldots, w_d)$ with the $d\times d$ matrices $(\partial
\rho/\partial w)$ and $(\partial
\rho/\partial t)$ invertible at $g(p)$, where $t=\Im w$.  Hence, by
the implicit function theorem, $g(U)$ is given near $g(p)$ by
\eqref{phi'} with $\phi =g_2\circ g_1^{-1}$.  We observe  that
since $\rho= t -  \phi(x,y,s)$ is also a vector-valued defining function of
$g(U)$ near $g(p)$, the invertibility of the matrix  $(\id + i
{\partial
\phi/ \partial s})$ follows from that of $(\partial\rho/\partial w)$.
The proof of the lemma is complete.
\epf

\br  For a germ of a generic embedding $g$, as in the statement of
Lemma \ref{germ}, even if $g(U)$ is given by an equation of the
form \eqref{phi'} for some coordinates $(z,w)$ in $\C^N$, it does
not necessarily follow that the matrix $(\id + i {\partial
\phi/ \partial s})$ is invertible.  For example, for the generic
embedding $g:\R^4\to\C^3$ with
$$ g(x,y,s_1,s_2) = \big(x+iy, s_1 +i(x+s_2), s_2 +i(y-s_1)\big),
$$
$g(\R^4)$ is given by \eqref{phi'} with $\phi(x,y,s)= (x+s_2,
y-s_1)$, but the $2\times 2$ matrix
$(\id + i {\partial
\phi/ \partial s})$ is not invertible.\er

The following definition will be needed for Lemma \ref{split} below.
\bd
Let $(\6M,p)$ and $(\6M',p')$ be two germs of smooth submanifolds of
$\C^N$ and $k$ a positive integer.  We shall say that $(\6M,p)$ and
$(\6M',p')$ are {\em $k$-equivalent} if there exists a germ of a local
biholomorphism $H:(\C^N,p)\to (\C^N,p')$ such that if $x \mapsto Z(x)$ is a
local parametrization of $\6M$ defined in a neighborhood of $0$ in
$\R^{\dim \6M}$ with $Z(0) = 0$ and
$\rho'$ is a vector-valued defining function for $\6M'$, then
$$
\rho'(H(Z(x)),\1{H(Z(x))}) = O(|x|^{k+1}).
$$
\ed
It is clear from the definitions that if $(\6M,p)$ and $(\6M',p')$ are generic
submanifolds of $\C^N$ which are $k$-equivalent, then $(\6M',p')$ is
$(k-1)$-nondegenerate (resp. of strong type $k$) if and only if the same
holds for $(\6M,p)$.

For a positive integer  $k $,  we split the jet space
$J^k_{0,0}(\C^n\times\R^d,\R^d)$ into its harmonic and
nonharmonic (free from harmonic terms) parts, i.e.
\begin{equation}\Label{harmfull}
J^k_{0,0}(\C^n\times\R^d,\R^d) =(J^k)_h \oplus (J^k)_{nh},
\end{equation}
  where
\begin{equation}\Label{harm} (J^k)_h:=\{j^k_0 \psi :
\psi=\Re(\sum_{0\ne |\a|+|\g|\le k}
\psi_{\a\g}z^\a s^\g),\, \psi_{\a\g}\in\C^d,\, z\in\C^n,\,
s\in\R^d\}
\end{equation}  and
\begin{equation}\Label{harmless} (J^k)_{nh}:=\{j^k_0
\psi :
\psi=\Re(\sum_{{\a\ne 0,
\b\ne 0}\atop{ |\a|+|\b|+|\g| \le k}}
\psi_{\a\b\g}z^\a \bar z^\b s^\g),\,
\psi_{\a\b\g}\in\C^d,\, z\in\C^n,\, s\in\R^d\}.
\end{equation}

We now state the main result of this section, a splitting of the jet
spaces of generic embeddings, which will be used in Sections
\ref{finitely} and \ref{s:strong}.

\bl\Label{split} Let $(z,w)$ be linear coordinates in
$\C^N$ ($N=n+d, m=2n+d$), $p\in \R^m$, and let
$\6O\subset J^k_p:= J^k_p(\R^{m},\C^N)$ be the Zariski
open subset of
$k$-jets of the form $j^k_pg$, where $g:(\R^{m},p)\to \C^N$
is a germ at $p$ of a generic embedding
  whose image is of the form
\eqref{phi'} near $g(p)$ with $(\id + i {\partial \phi/ \partial s})$
invertible (as in Lemma \ref{germ}). Then there exists a
birational map
\begin{equation}\Label{crap}
\Psi=(\Psi_1,\Psi_2,\Psi_3,\Psi_4)\colon J^k_p
\to\C^N
\times J^k_{0,0}(\R^m,\R^m) \times (J^k)_{h}\times
(J^k)_{nh}
\end{equation} smooth on $\6O$, such that, if $y=j^k_pg\in
\6O$ and
$\Psi_4(y)=j^k_0\psi$, then the germ  at $g(p)$ of
the image of $g$ in $\C^N$  is $k$-equivalent to the germ at
$0$ of the submanifold $\Im w=\psi(z,\bar z,\Re w)$.
Moreover, $\Psi$ is a diffeomorphism between $\6O$ and  a
Zariski open subset in the target space in
\eqref{crap}.
\el

\bpf We use the natural identification $J^k_p\cong \C^N\times
J^k_{0,0}(\R^m, \C^N)$ and the  notation
$y=(y_\a)_{0\le|\a|\le k}\in J^k_p$ with $y_\a\in\R^{2N}$,
$\a\in\Z^{m}_+$,
$y_\a=(y_\a^j)_{1\le j\le 2N}$. We set
$$\Psi_1(y):=y_0\in\R^{2N}\cong\C^N,
\quad \Psi_2(y):=(y_\a^j)_{1\le|\a|\le k,1\le j\le m}
\in J^k_{0,0}(\R^m,\R^m).$$ Next, we complete the pair
$(\Psi_1,\Psi_2)$ to a birational map
\begin{equation}\Label{coffee} (\Psi_1,\Psi_2,\Phi)\colon
\6O\to
  \C^N\times J^k_{0,0}(\R^m,\R^m) \times
J^k_{0,0}(\C^n\times\R^d,\R^d),
\end{equation}  as a diffeomorphism onto its image, by constructing
$\Phi$ as follows. Let $(y_\a)=(\d^\a g(p))$ for
$1\le|\a|\le k$ and some germ $g=(g_1,g_2)\colon (\R^m,p)\to
(\R^{2n+d}\times \R^d,0)$ as in Lemma \ref{germ}. Recall that
$g_1$ is an invertible germ at $p$ and
$\phi:=g_2\circ g_1^{-1}$ so that for $U$ a sufficiently small
neighborhood of $p$ in $\R^{2n+d}$,
$g(U)$ is given near $g(p)$ by
\eqref{phi'}.  We then set $\Phi(y):=j^k_0 \phi$, with 
that is, $\Phi(y) =  (\d^\a \phi(0))_{1\le|\a|\le k}$.  It is clear that
$\Phi(y)$ depends only on $y=j^k_p g$ and not on the
representative $g$, and that
$\Phi(y)$ is rational in view of the chain rule. Moreover,  since
$j^k_p g_2 = j^k_p (\phi\circ g_1)$, the inverse of the map
\eqref{coffee} is a polynomial map. Thus the map
\eqref{coffee} is a birational diffeomorphism onto its image.
In fact, the image of \eqref{coffee} is the Zariski open set
given by  $(Z,\mu,\nu)$  with $Z\in\C^N$, $\mu$ invertible in
$J^k_{0,0}(\R^m,\R^m)$ and
$\nu=j^k_0\phi\in J^k_{0,0}(\C^n\times\R^d,\R^d)$ with
$(\id+i\phi_s(0))$ (with $s\in\R^d$) being an invertible
$d\times d$ matrix.

We now set $\Psi_3(y)$ to be the component  of
$\Phi(y)$ in $(J^k)_h$ according to the decomposition
\eqref{harmfull}.  We may assume that
$\Phi(y)=j^k_0\phi$  with $\phi$ being a polynomial. In order
to define
$\Psi_4(y)$, we need to eliminate the harmonic components in
$\phi$ by an appropriate biholomorphic change of coordinates
near
$0$ in $\C^N$.  We look for
new coordinates
$(z',w')\in\C^n\times\C^d$ given by
$z=z'$, $w=h(z',w')$ with  $h(0)=0$,
$h_{w'}(0)\ne 0$. In these coordinates $g(U)$ is given by
\begin{equation}\Label{rho} h(z',w')-\bar h(\bar z',\bar w')
=2i\phi \big(z',\bar z',\frac{h(z',w') +\bar h(\bar z',\bar
w')}{2}\big).
\end{equation}  We may consider $z',\bar z', w', \bar w'$ as
independent variables.  In order that $g(U)$ be given in the
$(z',w')$ coordinates by an equation of the form
\eqref{phi'}, where the defining function has no nonzero
harmonic terms, it is necessary and sufficient that the equality
\eqref{rho} holds identically when $\bar z'=0$ and
$\bar w'=w'$ (see e.g.\ \cite{BER}).  That is, we must have
\begin{equation}\Label{rho1} h(z',w')-\bar h(0,w')
\equiv 2i\phi \big(z',0,\frac{h(z',w') +\bar h(0, w')}{2}\big).
\end{equation}  To solve \eqref{rho1} for $h$ we first set
$z'=0$ and
$w'=s\in\R^d$ to obtain
\begin{equation}\Label{rho2}
\Im h(0,s)
\equiv \phi (0,0,\Re h(0,s)).
\end{equation}  It suffices to take
\begin{equation}\Label{better} h(0,s):=s+i\phi(0,0,s).
\end{equation}
  With this choice, $h(z',w')$ is uniquely determined from
\eqref{rho1} by an immediate application of the implicit
function theorem (since the  $d\times d$ matrix $(\id+i\phi_s(0))$
is assumed invertible). In the new coordinates $(z',w')$ we can
write a defining equation of
$g(M)$ in the form
$\Im w'=\psi(z',\bar z',\Re w')$, with $j^k_0 \psi$ in
$(J^k)_{nh}$ (as given by \eqref{harmless}). Conversely,
knowing
$j^k_0\psi$  and the harmonic part
$\Psi_3(y)$ of $j^k_0\phi$, we can recover $j^k_0 h$ from
\eqref{better} and \eqref{rho1} and hence  the entire jet
$j^k_0\phi$. We may now define
$\Psi_4(y):=j^k_0\psi$, completing the proof of the lemma.
\epf

\section {Jet spaces of finitely nondegenerate embeddings}
\Label{finitely}
We need to identify those jets of germs
$g: (M,p) \to X$ of generic embeddings whose images are finitely
nondegenerate at $g(p)$. For this it suffices to consider the case
$M=\R^m$, with $m=\dim M = 2n+d$ and $X=\C^N$.  We fix an
integer
$k\ge 2$, and, as in the
previous section, we let $J^k_p:= J^k_p(\R^{m},\C^N)$.
The following stratification  result will be an important ingredient in
the proofs of our main results.

\bp\Label{wp} Let $W^0_p\subset J^k_p$ be the set of all
$k$-jets at $p$  of germs of generic embeddings $g\colon (\R^m,p)\to
\C^N$ and $W'_p\subset W^0_p$ be the subset consisting of
all those
$k$-jets of generic embeddings  that are
$(k-1)$-degenerate. Then $W^0_p$ is a Zariski open subset of
$J^k_p$ and  $W'_p$ is a real-algebraic subset of $W^0_p$
admitting the stratification into real-analytic submanifolds
$$W'_p = \bigcup_{1\le e\le n} W'_{p,e},$$
where $W'_{p,e}$ is the set of all $k$-jets of generic embeddings
whose $(k-1)$-degeneracy at
$g(p)$ (as defined by \eqref{rp}) is precisely $e$. 
In addition, for each $e$, $1\le e\le n$, the closure of
${W'_{p,e}}$ in ${W^0_{p}}$ admits the
stratification
$$\overline{W'_{p,e}} = \bigcup _{e\le c\le n} W'_{p,c}.$$
Furthermore, for $1\le e\le n$, the (real) codimension of the real-analytic submanifold
$W'_{p,e}$
  in
$J^k_p$ is at least
\begin{equation}\Label{minuses} 2d\left({k+n-1\atop
k-1}\right) -2n- 3d+2,
\end{equation}
and this bound is achieved for $e=1$.
\ep

\bpf
Recall that, if $g$ is defined in an open neighborhood  $U
$ of $0$ in $\R^m$ and $g(U)$ is given by
\eqref{phi'} with
$j^k_0\phi$ in $(J^k)_{nh}$ (as defined by
\eqref{harmless}), then $g(U)$ is $(k-1)$-degenerate at
$0$ if and only if the rank of the matrix $A$ whose rows are
\begin{equation}\Label{ja}
\phi^j_{z,\bar z^{\a}}(0)=(\phi^j_{z_1,\bar
z^\a}(0),\ldots,\phi^j_{z_n,\bar z^{\a}}(0)) ,\quad 1\le j\le d,
\quad 1\le|\a|\le k-1,
\end{equation}
is less than or equal $n-1$ (see e.g.\ \cite{BER}).  
We should observe here that the rank of
the matrix $A$ is independent of the choice of complex coordinates in
$\C^n$, as well as of the representative $\phi$ of the jet $j_0^k\phi$.
Since the number of multiindices
$\a\in\Z^n_+$ with
$0\le |\a|\le k-1$ is $\left({k+n-1\atop k-1}\right)$, it follows
that the number of rows of $A$ is
$d\left({k+n-1\atop k-1}\right)-d$. We will need the
following lemma.

\bl\Label{wp1'} Let $(J^k)_{nh}$ be defined by
\eqref{harmless} and, for $0\le r\le n$,  let $\mathcal
A_r\subset (J^k)_{nh}$ be the (semialgebraic) subset defined by
\begin{equation}\Label{A'}
\mathcal A_r:=\{j^k_0\phi\in (J^k)_{nh} : \rk A = r\},
\end{equation} where  $A$ is the matrix associated to
$j^k_0\phi$ whose rows are given by \eqref{ja}.
  Then $\mathcal A_r$ is  a real-analytic submanifold of
$(J^k)_{nh}$ with
\begin{equation}\Label{codimar}
\codim\, \mathcal A_r =  2d(n-r)\left({k+n-1\atop k-1}\right)
-(n-r)(2d+d(n-r)+2r).
\end{equation}
Moreover, if $\overline{\6A}_r$ denotes the closure of $\6A_r$ in
$(J^k)_{nh}$, then we have the stratification
\begin{equation}\Label{clo}
\overline{\6A}_r =\bigcup_{0\le c\le r}\6A_c.
\end{equation}
\el

\bpf  We fix $j^k_0\3\phi\in\mathcal A_r$ and denote by
$A_0$ the corresponding matrix.  Since $\rk A_0=r$, after a
complex-linear change of variables in $\C^n$, we may assume
that the kernel of $A_0$ is spanned by $e^a$ with
$a=r+1,\ldots,n$, where $e^1,\ldots,e^n$ are the standard
basis vectors in
$\C^n$. Furthermore, we can choose $r$ linearly independent
rows
$\3\phi^{j_1}_{z,\bar
z^{\a_1}}(0),\ldots,\3\phi^{j_r}_{z,\bar z^{\a_r}}(0)$  of
$A_0$. Moreover, since $e^a$  are in the kernel of $A_0$ for
$a=r+1,\ldots,n$, we have $\3\phi^j_{z_l,\bar z_a} =
\overline{\3\phi^j_{z_a,\bar z_l}}=0$, for $1\le l\le n$. Hence,
in the above choice of $\alpha_l$, $1\le l\le r$, we must
have
$\a_l\ne \eps^b$ for any $b=r+1,\ldots,n$, where
$\eps^b:=(0,\ldots,1,\ldots,0)$ with $1$ at the $b$-th position.

Then, for $j^k_0\phi$ near $j^k_0\3\phi$ in $(J^k)_{nh}$, the system
of linear equations in the unknowns
$(v_1,\ldots,v_n)$,
\begin{equation}\Label{vees}
\sum_{q=1}^n \phi^{j_l}_{z_q,\bar z^{\a_l}}(0) v_q=0,
\quad 1\le l \le r,
\end{equation} has an $(n-r)$-dimensional space of solutions,
whose basis can be chosen to be of the form
\begin{equation}\Label{norm}
v^i=(v^i_1,\ldots,v^i_{r},0,\ldots,1,\ldots,0)\in\C^n
\end{equation}
  with
$1$ at the $(r+i)$th position, where $i=1,\ldots,n-r$. Note that
$v^i_q$, $1\le i\le n-r$, $1\le q\le r$, are uniquely determined
by
\eqref{vees} and depend only on the chosen rows of $A$,
the matrix associated to the jet $j^k_0\phi$.


  To prove that
$\6A_r$ is a manifold of codimension $K_r$, where $K_r$
is the number given by the right-hand side of \eqref{codimar},
we shall show that
$\mathcal A_r$ is given near $j^k_0\3\phi$ by the vanishing
of $K_r$ real-analytic functions, whose differentials are
independent at every point near $j^k_0\3\phi$.

We now observe that a jet $j^k_0\phi\in (J^k)_{nh}$ near
$j^k_0\3\phi$ belongs to
$\6A_r$ if and only if $Av^i=0$ for $1\le i\le n-r$, where $A$
is the  matrix corresponding to $j^k_0\phi$ as usual. The latter
condition, in view of \eqref{vees}, is equivalent to the system,
\begin{equation}\Label{star}
\sum_{q=1}^n \phi^j_{z_q,\bar z^\a}(0) v^i_q = 0,
\end{equation} for
$1\le j\le d, \, 1\le i\le n-r,\, 1\le |\a|\le k-1$, where $(j,\a)\ne
(j_l,\a_l)$ for
$1\le l\le r$. Note that the coefficients $v^i_q$, $1\le i\le n-r$, $1\le
q\le r$, depend only on the variables $\phi^{j_l}_{z_q,\bar
z^{\a_l}}(0)$, $1\le l \le r$, $1\le q \le n$.  Hence
\eqref{star} can be considered as a linear system of equations
in the jet variables $\phi^j_{z_q,\bar z^\a}(0)$ for $j,q,\a$ as
in \eqref{star}.

Recall that we have $\a_l\ne\eps^b$ for any $b=r+1,\ldots,n$,
hence the system \eqref{star} contains, in particular, each
equation
\begin{equation}\Label{star1}
\sum_q \phi^j_{z_q,\bar z_{a+r}}(0)  v^i_q = 0, \quad  1\le
a,i\le n-r, \, 1\le j\le d.
\end{equation} We now consider a new system of equations
obtained from \eqref{star} by replacing each equation in
\eqref{star1} by the new equation
\begin{equation}\Label{replace'}  (\bar v^a_1,\ldots, \bar
v^a_n)
\begin{pmatrix}
\phi^j_{z_1,\bar z_1}(0) & \ldots & \phi^j_{z_n,\bar z_1}(0)
\\
\vdots & \ddots & \vdots \\
\phi^j_{z_1,\bar z_{n}}(0) & \ldots & \phi^j_{z_n,\bar
z_n}(0)  \\
\end{pmatrix}
\begin{pmatrix} v^i_1 \\
\vdots \\ v^i_n \\
\end{pmatrix} =0
\end{equation}  for the given $a,i,j$. Note that the equations
in the new system  obtained in this fashion are linear
combinations of equations in
\eqref{star}. Moreover, it is easy to see from the normalization
\eqref{norm} that the two systems are actually equivalent.

Recall that, in view of \eqref{norm}, $v^i_q=1$ for $q=r+i$.
Hence each equation in \eqref{star} with fixed $j,i,\a$, is
linear  in the variable $\phi^j_{z_{i+r},\bar z^\a}(0)$ with
coefficient $1$. Moreover, the latter variable is an arbitrary
complex number if either $|\a|\ge 2$ or $\a=\eps^a$ with $1\le
a\le r$. These correspond precisely to the equations of the old
system \eqref{star} that were not replaced. On the other hand,
each of  the new equations \eqref{replace'} is linear in
$\phi^j_{z_{i+r},\bar z_{a+r}}(0)$ with coefficient $1$.
Moreover, for $i=a$, both the variable $\phi^j_{z_{i+r},\bar
z_{a+r}}(0)$ and the equation \eqref{replace'} are real,
whereas for $i\ne a$, the corresponding equation
\eqref{replace'} is conjugate of that obtained by exchanging
$i$ and $a$.
In total, we obtain
$K_r$ independent real equations, proving the first part of the
lemma.

To prove \eqref{clo}, observe first that the rank of a matrix is
semi-continuous and hence the left-hand side of \eqref{clo}
is contained in the  right-hand side. For the converse inclusion,
let $j^k_0\phi\in \6A_c$ for $0\le c\le r$. As  in the beginning
of the proof of this lemma, we may assume that the kernel of the
corresponding matrix $A$ is spanned by the standard basis
vectors $e^a$, $c+1\le a\le n$.
Then, for any sufficiently small real
$\eps \ne 0$,  it is easy to see that $j^k_0(\phi+\eps
\2\phi) \in \6A_r$ with
$\2\phi(z,\bar z, s):=(\sum_{l=1}^{r} z_l\bar z_l,
0,\ldots,0)\in\R^d$.
This proves \eqref{clo} completing the proof of the lemma.
\epf

We now return to the proof of Proposition~\ref{wp}. We shall
first show that we can cover the set $W^0_p$  by finitely
many open subsets $\6O_l\subset J^k_p$ such that for every
$l$, there exist linear coordinates $(z^l,w^l)\in\C^N$ with
$j^k_p g\in \6O_l$ if and only if $g\colon (\R^m,p)\to \C^N$ is a
germ of a generic embedding whose image is of the form
\eqref{phi'}. Indeed, we start with a fixed set
$(Z_1,\ldots,Z_N)$ of linear coordinates in $\C^N$ and take
finitely many  new sets of linear coordinates
$(z^l,w^l)\in
\C^n\times \C^d$, obtained by permutations of the $Z_j$'s and
multiplication of some of the $Z_j$'s by
$i=\sqrt{-1}$, and take $\6O_l$ to be the set of all $k$-jets of
germs of generic embeddings $g$  whose images can be
graphed as in
\eqref{phi'} with respect to the coordinates $(z^l,w^l)$. Then
it is easy to see that $\6O_l$ is Zariski open and
$\cup_l\6O_l=W^0_p$. The rest of the proof follows from
Lemmas~\ref{split} and \ref{wp1'}. Indeed, for each choice of linear
coordinates
$(z^l,w^l)$ as above, we may apply Lemma~\ref{split}  to
obtain the map $\Psi^l$ as in \eqref{crap}, which is a
diffeomorphism on $\6O_l$. Then it easily
follows that
$W'_{p,e}=   \cup_l (\Psi^l_4|_{\6O_l})^{-1}
(\6A_{n-e})$  for each $e$, $1\le e \le n$,  where $\6A_r$ is given
by \eqref{A'}. This completes the proof of
Proposition~\ref{wp}, in view of Lemma~\ref{wp1'}.
\epf

\br
Note that the proof of Lemma~\ref{wp1'} actually
shows that $\6A_r$ contains the open dense subset
$$\3{\6A}_r:=\{j^k_0\phi\in \6A_r: \exists (b_1,\ldots,b_d)\in
\R^d\setminus \{0\}, \, \rk (b_1\phi^1_{z,\bar z}(0)+\ldots +
b_d\phi^d_{z,\bar z}(0)) = r
\},$$
where $\phi^j_{z,\bar z}(0)$ denotes the corresponding
Hermitian $n\times n$ matrix.
Note also that in general, $\3{\6A}_r$ is a proper subset of
$\6A_r$  even for $k=2$. Indeed, for $z\in\C^3$, $s\in\R$,
consider the  function
$$\phi(z,\bar z,s):= ( |z_1|^2 - |z_2|^2,
2\Re(z_1\bar z_3 + z_2\bar z_3) ).$$
Then $j^2_0\phi\in \6A_2$ but $j^2_0\phi\not\in \3{\6A}_2$,
i.e.\ the joint kernel of the Hermitian matrices $\phi^1_{z,\bar
z}(0)$ and
$\phi^2_{z,\bar z}(0)$ is zero but any  (real) linear combination
of them has a nontrivial kernel. Such an example is closely
related to the so-called ``null-quadrics'', i.e.\ quadrics of the form $\Im
w=H(z,\bar z)$, $(z,w)\in\C^n\times \C^d$, where $H=(H_1,\ldots,H_d)$
is a vector-valued Hermitian form with trivial common kernel but such
that any real linear combination of the $H_l$'s is a degenerate Hermitian
form. See \cite{Bel} for further details.
\er

\section{Jet spaces of embeddings of finite strong type}\Label{s:strong}

In this section we fix an integer $k \ge 2$ and identify those $k$-jets of
germs
$g: (M,p) \to X$ of generic embeddings whose images are not of strong
type $l$ at
$g(p)$ for any $l\le k$. In this context we shall prove a stratification
similar to that given by Proposition \ref{wp} (see Proposition \ref{ft}
below).  We start with the following lemma.

\bl\Label{type} Assume that $(\6M,0)\subset \C^N$ is a germ of
a generic submanifold given by
\eqref{phi'} near the origin with $j^k_0\phi\in (J^k)_{nh}$
(as defined by \eqref{harmless}). Then the codimension of the
span of  the vectors  \eqref{comm}  and their conjugates at
$p=0$ in $\C T_0 \6M$ coincides with that of the span of all
vectors
$\phi_{z_r, \bar  z^\a}(0),  \phi_{\bar z_r, z^\a}(0)$ in
$\C^d$  for $1\le r\le n$, $1\le |\a|\le k-1$.
\el

\bpf If $\6M$ is given by \eqref{phi'}, a local basis of the
$(1,0)$ vector fields can be chosen to be
$$L_j=\frac{\d}{\d z_j} + i\phi_{z_j}(z,\bar z,s)
(\id-i\phi_s(z,\bar z,s))^{-1}
\frac{\d}{\d s},
\quad j=1,\ldots,n.$$ Here we view $\phi_{z_j}$ as a row
vector in $\C^d$,
$\phi_s$ as a $d\times d$ matrix with $s=\Re w\in\R^d$ and
$\frac{\d}{\d s}$ as a column vector with $d$ components.
We shall prove, for $j_1,\ldots,j_r,l\in\{1,\ldots,n\}$ and
$r\le k-1$, the following identity, which will imply the
conclusion of the lemma:
\begin{equation}\Label{ident}
[L_{j_1},[L_{j_2},\ldots,[L_{j_r},\1L_l]\ldots]](0)=
-2i\phi_{z_{j_1},z_{j_2},\ldots,z_{j_r},\bar
z_l}(0)\frac{\d}{\d s}.
\end{equation} Indeed,  using  induction on $r$, it follows
from the assumption
$j^k_0 \phi\in(J^k)_{nh}$, that
\begin{equation}\Label{ident1}
[L_{j_1},[L_{j_2},\ldots,[L_{j_r},\1L_l]\ldots]](z,\bar z,s)=
\big(-2i\phi_{z_{j_1},z_{j_2},\ldots,z_{j_r},\bar z_l}(z,\bar
z,s) + \langle \bar z \rangle\ + O(k-r)
\big)\frac{\d}{\d s} ,
\end{equation} where $\langle \bar z\rangle$ denotes a
$\C^d$-valued row of polynomials of the form $\sum\bar z_j
f_j(z,\bar z,s)$ and $O(k-r)$ a vector-valued function
vanishing at the origin of order at least   $k-r$. This proves the
lemma.
\epf

\bp\Label{ft} As in Proposition \ref{wp}, let $W^0_p\subset J^k_p$
be the set of all
$k$-jets at $p$  of germs of generic embeddings $g\colon (\R^m,p)\to
\C^N$.  Denote by $W''_p\subset W^0_p$  the subset consisting
of all $k$-jets of generic embeddings  that are not of strong type
$l$ for any $l\le k$. Then $W''_p$ is a real-algebraic subset of
$W^0_p$ admitting the stratification into real-analytic submanifolds
$$W''_p = \bigcup_{1\le e\le d} W''_{p,e},$$
where $W''_{p,e}$ is the set of all $k$-jets of generic embeddings
whose $k$-defect (as defined by \eqref{defect}) is precisely $e$ at
$g(p)$. Furthermore the closure of $W''_{p,e}$ in $J^k_p$ admits the
stratification
$$\overline{W''_{p,e}} = \bigcup _{e\le c\le d} W''_{p,c}$$
and the (real) codimension of $W''_{p,e}$ for $1\le e\le d$, in
$J^k_p$ is at least
\begin{equation}\Label{pluses}
\Big[ 2n\Big({k+n-1\atop k-1}\Big) -n^2-2n- d+1\Big]^+
\end{equation} with the notation $[r]^+:=\max(r,0)$ for
$r\in\Z$,
and this bound is achieved for $e=1$.
\ep

\bpf If $g$ is defined in open neighborhood $U$ of $\R^m$, and
$g(U)$ is given by
\eqref{phi'} with
$\phi$ in
$(J^k)_{nh}$, then, by Lemma~\ref{type},
$g(U)$ is not of strong type $k$ at
$0$ if and only if the rank of the matrix $B$ whose rows are
\begin{equation}\Label{ja1}
\phi_{z^\a,\bar z_r}(0)=(\phi^1_{z^\a,\bar
z_r}(0),\ldots,\phi^d_{z^\a,\bar z_r}(0)) ,\quad 1\le r\le n,
\quad 1\le|\a|\le k-1,
\end{equation} and their conjugates,  is less than $d$.   Since
the number of multiindices $\a\in\Z^n_+$ with
$0\le |\a|\le k-1$ is $\left({k+n-1\atop k-1}\right)$, it follows
that the number of rows of $B$ with $|\a|>1$ is
$2n\left(\left({k+n-1\atop k-1}\right)-n-1\right)$. Since each
row appears together with its conjugate,  we don't change the
rank of $B$ by replacing the pair of each row and its conjugate
by the real and imaginary parts of that row. Hence we obtain
$2n\left(\left({k+n-1\atop k-1}\right)-n-1\right)$  real rows
for $|\a|>1$. Furthermore, each row in \eqref{ja1} with
$|\a|=1$ appears twice in the matrix $B$. In fact, the rows
$(\phi^1_{z_r,\bar z_r}(0),\ldots,\phi^d_{z_r,\bar z_r}(0))$, $1\le r\le
n$, are real and hence coincide with their conjugates, whereas the
conjugate of each row
$(\phi^1_{z_r,\bar z_l}(0),\ldots,\phi^d_{z_r,\bar z_l}(0))$
with $r\ne l$, appears again as the row
$(\phi^1_{z_l,\bar z_r}(0),\ldots,\phi^d_{z_l,\bar z_r}(0))$.
Hence the rank of $B$ does not change if we  remove these
$n^2$ repeated rows. We thus end up with a real matrix $B'$
having the same rank as $B$ with
$K'':=2n\left(\left({k+n-1\atop k-1}\right)-n-1\right)+n^2$
real rows, whose entries can be seen as part of independent
coordinates of $j^k_0\phi$ in $(J^k)_{nh}$.

We now proceed as in the proof of Proposition~\ref{wp} and
Lemma~\ref{wp1'}. We make use of the following lemma.

\bl\Label{wp2'}
Let $(J^k)_{nh}$ be defined by
\eqref{harmless} and, for $0\le r\le d$,  let $\mathcal
B_r\subset (J^k)_{nh}$ be the (semialgebraic) subset defined by
\begin{equation}\Label{A''}
\mathcal B_r:=\{j^k_0\phi\in (J^k)_{nh} : \rk B = r\},
\end{equation}
where  $B$ is the matrix associated to
$j^k_0\phi$ whose rows are given by \eqref{ja1}.
  Then $\mathcal B_r$ is  a real-analytic submanifold of
$(J^k)_{nh}$ with
\begin{equation}\Label{plusesb}
\codim\, \6B_r = (d-r)\Big[ 2n\Big({k+n-1\atop k-1}\Big) -n^2-2n-r
\Big]^+
\end{equation} with the notation $[r]^+:=\max(r,0)$ for
$r\in\Z$.
Moreover, if $\overline{\6B}_r$ denotes the closure of $\6B_r$ in
$(J^k)_{nh}$, then we have the stratification
\begin{equation}\Label{clo'}
\overline{\6B}_r =\bigcup_{0\le c\le r}\6B_c.
\end{equation}
\el

\bpf
It follows from the discussion preceeding the lemma that we can replace
the matrix $B$ by the real matrix $B'$ in \eqref{A''}.
Then each $\6B_r$ can be regarded as an orbit in the space of all  real
$K''\times d$ matrices under the action of ${\bf GL}(K'',\R)\times
{\bf GL}(d,\R)$ by right and left multiplication.
Hence each $\6B_r$ is a manifold of codimension $(d-r)[K''-r]^+$
which yields \eqref{plusesb}. We leave the remaining details of the proof
to the reader.
\epf

The rest of the proof of the proof of
Proposition~\ref{ft} follows closely that of Proposition~\ref{wp}.
\epf

\section{Proof of Theorems~\ref{weak} and
\ref{weak-precise}} \Label{proofweak}

\bpf
We start with the proof in the case where $M$ is a smooth manifold
and
$\tau: M \to X$ is a smooth map. We shall make use of Thom
transversality  theorem (see e.g.\ \cite{GG} or  \cite{Hir}).
We first choose $k_1=k_1(m,N)$  as defined by \eqref{k1} and fix $k:=k_1+1$.
We shall take
$W'\subset J^k(M,X)$ to be the subset of all
$k$-jets of the form $j^k_pg$, $p\in M$, with $g\colon
(M,p)\to X$ a germ of a generic embedding whose image is
($k-1$)-degenerate at
$g(p)$.  Observe that  the latter condition imposed on $g$
depends only on $j^k_p g$ and not on the representative $g$ of
this jet.

We now introduce a natural stratification of the set $W'$.
For each integer $1\le c\le n$, let $W'(c)\subset J^k(M,X)$
  be the set of all $k$-jets of germs of generic embeddings
$g\colon (M,p)\to X$ whose image has $(k-1)$-degeneracy at $g(p)$
equal to $c$.
Note that
\begin{equation}\Label{W}
W'=\bigcup_{1\le c\le n} W'(c).
\end{equation}
We claim that each $W'(c)$ is a smooth submanifold of
$J^k(M,X)$ of codimension greater than or equal to the integer given
by   \eqref{minuses}.    Moreover,
the  bound given by \eqref{minuses}  is achieved
for $c=1$.
Indeed, since the $(k-1)$-degeneracy is invariant under translations in the
target  space
$\C^N$, in the notation of Proposition~\ref{wp} we have
$$W'_{p,c}=W'_{(p,0),c}\times \C^N \subset J^k_{p,0}(M,\C^N)
\times \C^N,$$
with $W'_{(p,0),c}\subset J^k_{p,0}(M,\C^N)$
being a real-analytic submanifold of the same codimension as $W'_{p,c}$
in $J^k_{p}(M,\C^N)$.
If $U$ and $V$ are local real and complex coordinate
charts on $M$ and $X$ respectively, then  it is standard (see
e.g.\ \cite{GG}) that
$J^k(U,V)$ can be naturally identified with
$U\times V\times J^k_{0,0}(\R^m,\C^N)$, where $m:=\dim M=2n+d$ and
$N:=\dim_\C X=n+d$.
Then the claim for $W'(c)$ follows from Proposition~\ref{wp}.

For $g\in C^\infty(M,X)$, let $j^kg:=\{(p,j^k_pg)\in J^k(M,X): p\in M\}$.
Recall that it follows from Thom's transversality theorem that, if $M$ and
$X$ are as in Theorem~\ref{weak} and
$W\subset J^k(M,X)$ a real smooth submanifold whose codimension is greater
than
$\dim M$, then
$$\{g\in C^\infty(M,X) : j^k g\cap  W=\emptyset \}$$
  is  a
{\em residual} subset in $C^\infty(M,X)$, i.e.\ a countable intersection of
open dense subsets.
Here $C^\infty(M,X)$ is equipped with the Whitney topology.
For more details, see e.g.\ \cite[Chapter  2, Theorem 1.2]{Hir},
\cite[Chapter  2, Theorem 4.9]{GG}.

We now apply Thom's transversality
theorem to each submanifold $W'(c)$,
$1\le c\le n$, which are the strata of $W'$
given by \eqref{W}.
With the choice $k= k_1+1$, it follows that $j^k_p g\in W'$ if and only if
$g(M)$ is $k_1$-degenerate at
$g(p)$. One can check that if $k_1$ is the integer given by \eqref{k1} and
$k=k_1+1$, then the integer given by
\eqref{minuses} is strictly greater than $\dim M = 2n+d$.  Hence, by Thom's
transversality theorem, the set
$$\{g\in C^\infty(M,X) : j^k g\cap  W'=\emptyset \}$$ is
residual in $C^\infty(M,X)$.

To approximate by embeddings whose images are of strong type, we now let 
$k: = k_2(m,N)$,
  the integer given by \eqref{k2}.  Observe that for this choice, the
number given by
\eqref{pluses} is greater than $\dim M=2n+d$.  We let $W''\subset J^k(M,X)$
be the subset of all jets $j^k_pg$ with $g:(M,p) \to X$ a germ of a generic
embedding whose image is not of strong type $\le k_2$ at $g(p)$. As we did
for $W'$, we may stratify $W''$ and write
\begin{equation}
W''=\bigcup_{1\le e\le d} W''(e),
\end{equation}
  where $W''(e) \subset J^k(M,X)$ is the set of all $k$-jets of germs of generic
embeddings $g:(M,p)\to X$ whose image has $k$-defect at $g(p)$ equal to
$e$. It easily follows from Proposition \ref{ft} that each stratum $W''(e)$ is
a smooth submanifold of $J^k(M,X)$ whose codimension is at least the
number given by \eqref{pluses} and that this bound is achieved for $e=1$. We
again apply Thom's transversality  theorem to each stratum
$W''(e)$ of
$W''$ to conclude that
$$\{g\in C^\infty(M,X) : j^k g\cap  W''=\emptyset \}$$ is
also residual in $C^\infty(M,X)$.

Since the intersection of two residual sets is again residual and, in
particular, dense, and since the generic embeddings from $M$ to $X$ form
an open subset in $C^\infty(M,X)$
we obtain the approximation property as stated in
Theorem~\ref{weak-precise} for the smooth case. Since strong type $\ell$ for
some positive
$\ell$ implies finite type, we have also proved Theorem~\ref{weak} for the
smooth case.

For the real-analytic case, i.e.\ when both $M$ and $\tau$ are
real-analytic, let $U\subset C^\infty(M,X)$ be any open neighborhood of
$\tau$ (in the Whitney topology).  Note that the set of all generic embeddings
$\sigma$ in $U$ that are $\ell_1$-nondegenerate, for some $\ell_1\le k_1$,
and of strong type $\ell_2$, for some $\ell_2\le k_2$, is an open
subset of
$U$ in the Whitney topology, which is nonempty, by the proof above for the
smooth case. Since the set of all real-analytic maps from $M$ to $X$ is  dense
in
$C^\infty(M,X)$ (see e.g.\ \cite{Hir}), the desired generic embedding $\sigma$
can also be chosen to be real-analytic.
This completes the proofs of Theorems~\ref{weak} and
~\ref{weak-precise}.

\epf

\section{Proof of Theorems~\ref{strong} and
\ref{strong-precise}}\Label {proofstrong}
\begin{proof}The proof will follow in several steps.  In Step 1, we
consider the case
  $X=\C^N$ and
$M$ an open set in $\R^m$.  In Step 2, we still take $X=\C^N$,
but allow $M$ to be a general abstract manifold.  Step~3 deals with the
general case.

\subsection{Step 1}
Assume first that
$X=\C^N$ and
$M$ is an open set in $\R^m$ containing the origin. For any integer $k> 1$, we
use the  shorthand notation and the identifications
$J^k_p=J^k_p(M,\C^N)\cong \ J^k_{0}(\R^m,\C^N)$
and $J^k=J^k(M,\C^N)\cong M\times
J^k_{0}(\R^m,\C^N)$.
We
look for $\2\tau(p,t)$ of the form
\begin{equation}\Label{look}
\2\tau(p,t):=
\tau(p)+t  f(p), \quad t\in (-1,1),\, f\in C^\infty(M,\C^N).
\end{equation}
  As
in the proof of Theorem~\ref{weak}, we shall use  Thom's
transversality theorem to select $f$ appropriately.
For this, we first choose
$k=k_1+1$, with $k_1 = k_1(m,N)$ given by \eqref{k1}, and let $W'_p \subset
J^k_p$ be the semialgebraic subset defined in Proposition~\ref{wp}.
We shall show that the set
\begin{equation}\Label{res1}\{f\in C^\infty(M,\C^N): j_p^k\tau + tj^k_pf
\not\in W_p', \
\forall p \in M,\ \forall  t \in (-1,1)\setminus\{0\}\}
\text{ with } k=k_1+1,
\end{equation}
is residual.

Similarly, we  take $k = k'_2$, where $k'_2$ is the smallest integer
$k$ for which inequality \eqref{k2'} holds, and let
$W''_p\subset J^k_p$ be the semialgebraic subset defined in
Proposition~\ref{ft}. We shall again show that the set
\begin{equation}\Label{res2}\{f\in C^\infty(M,\C^N): j_p^k\tau + tj^k_pf
\not\in W_p'',
\
\forall p \in M,\ \forall  t \in (-1,1)\setminus\{0\}\} \text{ with } k=k'_2,
\end{equation}
is also residual.

We shall now complete the proof of  Theorem \ref{strong-precise} in the case
where $M$ is an open subset of $\R^m$ and $X=\C^N$. As in the statement of
the theorem, let $U$ be an open neighborhood of $\tau$ in $C^\infty(M,X)$.
Since $\tau$ is a generic embedding, and the set of all generic embeddings of
$M$ into $X$ is open  in $C^\infty(M,X)$, we may assume, by shrinking $U$,
if necessary, that all elements in $U$ are also generic embeddings. It follows
from the definition of the Whitney topology that there exists a continuous
positive function
$\delta(p)$ on $M$ and a positive integer  $r$ such that the set
$$U_\delta^r:= \{ g\in C^\infty(M,\C^N): |j_p^rg - j_p^r \tau|< \delta(p),\
\forall p
\in M\}
$$
is contained in $U$.  Here $|\cdot|$ is any norm on $J_0^r(M,\C^N)$, with the
identification given above.  Choose $f \in C^\infty(M,\C^N)$ such that
$\tau + f
\in U^r_\delta$ and $f$ in each of the residual sets defined by 
\eqref{res1} and
\eqref{res2}, which is possible since the intersection of two residual sets is
residual, and hence dense.  It follows from the definition of $U_\delta^r$ that
$\tau+tf$ is also in $U^r_\delta \subset U$ for any  $t \in [-1,1]$. 
It is now clear
that $\2\tau(p,t) = \tau(p) + tf(p)$, with $f$ as chosen above, satisfies the
conclusion of Theorem \ref{strong-precise}.

To complete Step 1 of the proof, it remains to show that the sets
\eqref{res1} and \eqref{res2} are residual.
  Note that with the identification above, we have
\begin{equation}\Label{w-equal}
W'_p=W'_0, \quad W''_p=W''_0, \quad p\in M.
\end{equation}
Recall that a semialgebraic set $A$ admits a finite (semialgebraic)
stratification into a disjoint union of real-analytic submanifolds, and the
maximum stratum dimension (resp.\ minimum stratum codimension) is
independent of the stratification and is said to be the dimension
(resp.\ codimension)  of $A$ (see e.g\ \cite{BR}).  We consider first the
  set given by
\eqref{res1}, for which $k: = k_1 + 1$.

  For $p\in
M$, let
$S_p\subset (J^k_p)^2\times ((-1,1)\setminus\{0\})$ be given by
$$S_p:=\{(\L_0,\L,t) : \L_0,\L\in J^k_p, \, t\in (-1,1)\setminus\{0\}, \,
\L_0+t\L\in W'_p \big\}.$$
By Proposition~\ref{wp},
for any $p$,
$S_p$ is semialgebraic with  $\codim\, S_p=\codim\, W'_p$.
Consider the following natural projections
$$\pi_p\colon (J^k_p)^2\times ((-1,1)\setminus\{0\}) \to (J^k_p)^2,
\quad \rho_p\colon  (J^k_p)^2\times ((-1,1)\setminus\{0\}) \to
J^k_p,
\quad \sigma_p\colon (J^k_p)^2 \to J^k_p$$
given by
$$\pi_p(\L_0,\L,t) = (\L_0,\L),
\quad \rho_p(\L_0,\L,t) = \L_0,
\quad \sigma_p(\L_0,\L) =\L_0.$$
Then we have
$$\codim_{\rho_p^{-1}(\L_0)} (\rho_p^{-1}(\L_0) \cap S_p)
= \codim_{J^k_p} W'_p.$$
Furthermore, $B_p:=\pi_p(S_p)\subset (J^k_p)^2$
   is also
semialgebraic by a theorem of Tarski-Seidenberg (see  e.g.\ \cite{BR}).
Since the dimension of a semialgebraic set cannot increase after projection
(see  e.g.\
\cite{BR}), we have
$$\codim\, B_p\ge \codim \, S_p-1=\codim \,W'_p-1.$$
We also have for every $\L_0\in J^k_p$,
$$ \pi_p(\rho_p^{-1}(\L_0) \cap S_p) = \sigma_p^{-1}(\L_0)\cap B_p, $$
implying
$$\dim (\sigma_p^{-1}(\L_0)\cap B_p)\le \dim (\rho_p^{-1}(\L_0) \cap
S_p),$$
and hence
\begin{equation}\Label{B}
\codim_{\sigma_p^{-1}(\L_0)} (\sigma_p^{-1}(\L_0)\cap B_p)\ge
\codim_{\rho_p^{-1}(\L_0)} (\rho_p^{-1}(\L_0)
\cap S_p) -1 = \codim \,W'_p-1.
\end{equation}

It follows from \eqref{w-equal} that $S_p=S_0$ and $B_p=B_0$.
Since $B_0$ is a semialgebraic set, we may consider its stratification into a
finite union of real-analytic disjoint submanifolds $B_0=\cup_j B_0^j $.
For clarity and motivation of the proof, we begin with the
simplifying assumption that $B_0$ consists of a single stratum, i.e.\ $B_0$
is a real-analytic submanifold of $(J^k_0)^2$, and that
$\sigma_p|_{B_0}$ is a submersion onto an open subset
$A_0\subset J^k_0$ (i.e.\ $\sigma_p|_{B_0}$ is of maximal
rank equal to
$\dim_\R J^k_0$ at every point of $B_0$).
We set
$$B:= \cup_{p\in M} B_p = M\times B_0 \subset \cup_{p\in
M}(J_p^k)^2,\quad
A:= M\times A_0 \subset J^k.$$
Then the map
$$\sigma \colon M\times (J^k_0)^2 \to M\times J^k_0 = J^k, \quad
\sigma(p,\L_0,\L):=(p,\L_0),$$ when restricted to
$B$, is a submersion onto the open subset
$A\subset J^k$. Define the map
$\l_0\colon M\to J^k_0$ by $\l_0(p):= j^k_p \tau \in J^k_0\cong J^k_p$,
where
$\tau\colon M\to \C^N$ is the given generic embedding.

Recall that $k_1$ is chosen to be $k-1$, where $k$
is the minimum integer
for which
$\codim W'_p$ (given by \eqref{minuses}) is greater than
$\dim_\R M (=2n+d)$.
The reader can also check that the same $k$ is also the minimum integer,
for which
$\codim W'_p -1 > \dim_\R M$.

We now consider the subset $V\subset J^k\cong M\times J^k_0$ defined
by
\begin{equation}\Label{V} V:=\{(p,\L)\in J^k :
(\l_0(p),\L)\in B_0 \},
\end{equation}
which is a submanifold of $J^k$ of
codimension
$\codim_{\sigma_0^{-1}(\L_0)} B_0$.
Hence by \eqref{B} and the choice of $k$, we have
$\codim V \ge \codim W_p -1 > \dim_\R M$. Thus we can apply
Thom's transversality theorem to   $V$ to obtain a residual set of
functions $f\in C^\infty(M,\C^N)$ with
$j^k f\cap V=\emptyset$, i.e.\
$j^k_p\tau + t j^k_p  f\notin W'_p$ for all $t\in (-1,1)\setminus \{0\}$
  and all $p\in M$.
This completes
the proof in the simplifying case where
$B_0\subset (J^k_0)^2$ is a real-analytic submanifold and
$\sigma_0|_{B_0}$ is a submersion onto an open subset
$A_0\subset J^k_0$.

We now consider the less restrictive assumption that
$B_0$ is still a real-analytic submanifold of $J^k_0$
(i.e.\ a single stratum)
  and $\sigma_0|_{B_0}$ is a
submersion  onto a submanifold
$A_0\subset J^k_0$  which is not necessarily open.
It should be noted that in
this case, the set $V$ defined by \eqref{V} may not necessarily be a
submanifold of $J^k$.
To remedy this,  we will enlarge  both
$A_0$ and $B_0$  as follows.
  Let $\Omega_0\subset J^k_0$  be an open neighborhood of
the submanifold $A_0\subset J^k_0$ with a retraction $r_0\colon
\Omega_0\to A_0$, i.e.\
$r_0\in C^\infty(\Omega_0,A_0)$ and
$r_0(\L_0)=\L_0$ for $\L_0\in A_0$.
The existence of such $\Omega_0$ and
$r_0$ is well-known (see e.g.\ \cite[Chapter 2, \S7]{GG}).
  We now define the enlargement of $B_0$, denoted by
$\2B_0\subset (J^k_0)^2$, to be
$$\2B_0:=\{(\L_0,\L) : \L_0\in \Omega_0, \, (r_0(\L_0),\L)\in B_0\} .
$$ Then it is easy to check that $\2B_0$  is a submanifold of $(J^k_0)^2$
containing $B_0$ and $\sigma_0|_{\2B_0}$ is a submersion onto the
open subset $\Omega_0$.
Moreover, it is also easy to check that
$$\codim_{\sigma_0^{-1}(\L_0)}(\sigma_0^{-1}(\L_0) \cap \2B_0) =
\codim_{\sigma_0^{-1}(\L'_0)}(\sigma_0^{-1}(\L'_0) \cap B_0)
\ge \codim W_0 -1$$
for any $\L_0\in \Omega_0$ and $\L'_0\in A_0$,
where the latter inequality follows from \eqref{B}.
(The reader should note that if $A_0$ is not open in $J^k_0$, then
$\dim \2B_0 > \dim B_0$.)
The rest of the proof
in this case can be reduced to the previous case above
(with $A_0$ open in $J^k_0$) by replacing
$B_0$ with $\2B_0$.

Finally, we consider the general case, i.e.\ when no restrictions are
imposed on $B_0$; it is merely semialgebraic and hence is a finite union
of real-analytic strata.  Then the proof is reduced to the previously
considered two cases by applying the following general result about
semialgebraic sets that follows by induction  from the
stratification of semialgebraic sets and the theorem of Tarski-Seidenberg
mentioned above.

\bl\Label{strat}
Let $C\subset \R^{K_1}\times \R^{K_2}$
be a semialgebraic subset
and $\sigma\colon \R^{K_1}\times \R^{K_2}\to \R^{K_1}$
be the canonical projection. Then
$C$ can be decomposed into a finite union of disjoint
real-analytic submanifolds $C_j$ such the restriction of
$\sigma$ to each $C_j$ is a submersion onto a real-analytic
submanifold of $\R^{K_1}$.
\el

\bpf[Sketch of the proof of Lemma \ref{strat}]
The proof is by induction on the dimension of $C$.
Since $C$ can be stratified as a finite union of real-analytic connected
submanifolds, we may take a stratum $C_0$ of maximal dimension and
assume, after removing a semialgebraic subset of lower dimension,
that this stratum is globally defined by the vanishing of finitely many
polynomials with independent differentials.
Furthermore, after removing another lower-dimensional semialgebraic
subset, we may assume that the restriction of $\sigma$ to $C_0$ is of
constant rank. By a theorem of Tarski-Seidenberg, $\sigma(C_0)$ is again
semialgebraic and hence can be stratified. Let $D_0\subset \sigma(C_0)$
be the union of strata of maximum dimension.
Then $\2C_0: =\sigma^{-1}(D_0)$ is an open submanifold of $C_0$,
whose complement is semialgebraic of lower dimension and such
that the restriction of $\sigma$ to $\2C_0$ is a submersion onto a
real-analytic submanifold.
\epf

The proof that the set given by \eqref{res1} is residual is complete
under the hypotheses of Step 1.  The proof that the set given by
\eqref{res2} is also residual follows by repeating verbatim the proof for
\eqref{res1} by replacing $W'_p$ by $W''_p$ and using Proposition
\ref{ft} in place of Proposition \ref{wp}.
The proofs of Theorems~\ref{strong} and ~\ref{strong-precise}
in the case when $M$ is an open set in $\R^m$ and $X=\C^N$ are now
complete.

\subsection{Step 2}
We now consider the case when $M$ is a general smooth manifold but
still assume $X=\C^N$.
We again look for a suitable deformation $\2\tau$ of $\tau$,  of the form
\eqref{look}, where we choose $f\in C^\infty(M,\C^N)$ appropriately. As
in  Step 1, we must show that the sets given by \eqref{res1} and
\eqref{res2} are residual.  For this, we cover
$M$ by countably many coordinate charts $M_j$.
Then the arguments of
Step 1, applied first to $k=k_1+1$ (as in Step 1), to each $M_j$,
yield subsets
$V_j\subset J^k(M_j,\C^N)$ such that
each $V_j$ is a finite union of smooth manifolds and, for $f\in
C^\infty(M,\C^N)$,
$(p,j^k_p f)\notin V_j$ for $p\in M_j$ implies that $\2\tau(\cdot,t)$ (given
by \eqref{look}),\ restricted to $M_j$ is a generic embedding, whose image
is
$k_1$-nondegenerate.
Then we regard each $V_j$ as a finite union of smooth submanifolds of
$J^k(M,\C^N)$ and apply Thom's transversality theorem to the union of
the
$V_j$'s in $J^k(M,\C^N)$ to obtain that \eqref{res1} is residual. An
identical argument with $k = k_2'$ shows that \eqref{res2} is residual.
The proof is then completed as in Step 1.  At this point we have proved
Theorems \ref{strong} and \ref{strong-precise} in the case where $X =
\C^N$, including the estimates for $k'_1$ and $k'_2$.

\subsection{Step 3}
We now treat the general case where $X$ is a complex manifold
and we are given a neighborhood $U$ of $\tau$ in $C^\infty (M,X)$.
By possibly shrinking $U$, we may assume that all maps in $U$ are
generic embeddings of $M$ in $X$.
We begin with locally finite coverings
of
$M$ and
$X$  by countably many
relatively compact coordinate
charts
$$M=\bigcup_j M_j= \bigcup_j \2M_j,
\quad X=\bigcup_j X_j, \quad j = 1, 2,\ldots,$$
with $M_j\Subset \2M_j\Subset M$ and $\tau(\2M_j)\Subset X_j\Subset
X$.
We now choose $k$ such that both integers \eqref{minuses} and
\eqref{pluses} are greater than $\dim M+2$.
For every $p\in M$, we consider the semialgebraic sets $W_p\subset
J^k_p(M,X)$ consisting of $k$-jets of germs of smooth generic
embeddings $g\colon (M,p)\to X$ whose images are either
$(k-1)$-degenerate or not of strong type $\le k$.
With this choice, we have
\begin{equation}\Label{codimw}
\codim W_p>\dim_\R M+2
\end{equation}
  as consequence of Propositions~\ref{wp} and
~\ref{ft}.

We also set $M_0=\2M_0 :=\emptyset$, $X_0:=\emptyset$ and shall
construct inductively a sequence of deformations
$\2\tau_j\colon M\times (-1,1)\to X$, $j=0,1,\ldots$, with the following
properties:
\begin{enumerate}
\item $\2\tau_j(p,0)\equiv \tau(p)$ for all $p\in M$;
\item $j^k_p\2\tau_j(\cdot,t)\notin W_p$ for $p\in
\bigcup_{l\le j} \1M_l$ and $t\in(-1,1)\setminus \{0\}$;
\item $\2\tau_j (\2M_l\times (-1,1))\Subset  X_l$ for all $l$;
\item $\2\tau_j(\cdot,t)\in U$ for all $t\in(-1,1)$;
\item $\2\tau_j(p,t)\equiv \2\tau_{j-1}(p,t)$ for $p\notin \2M_j$, $t\in
(-1,1)$, and $j \ge 1$.
\end{enumerate}
Then it is clear that $\2\tau_0(p,t):=\tau(p)$ satisfies the hypotheses (1)--(5)
for $j=0$.

We now assume that the $\2\tau_j$ have been chosen for $0\le j\le j_0$,
and construct $\2\tau_{j_0+1}$. For this, we slightly modify the
end of the proof in Step 1 by introducing an additional real parameter.
In that proof we take $M$ to be $\2M_{j_0+1}$,
regarded as a open set in $\R^m$,
and choose the same sets $A$, $B$, $A_0$, $B_0$.
We first look for $\3\tau$ of the form
\begin{equation}\Label{nouse}
\3\tau(p,t,t')= \2\tau_{j_0}(p,t') + t f(p,t'), \quad p\in \2M_{j_0+1},
\end{equation}
where $f\in C^{\infty}(\2M_{j_0+1}\times (-1,1), \C^N)$, with the
addition and multiplication by
$t$ understood to be with respect to
the chosen coordinates in $X_{j_0+1}$.
This is possible by the assumption (3) for $\2\tau_{j_0}$.

We now define $\l_0\colon \2M_{j_0+1}\times(-1,1)\to J^k_0$ by
$\l_0(p,t'):= j^k_p \2\tau_{j_0}(\cdot,t')$ and follow the arguments of Step
1 as follows. We fix a finite stratification $(B^l_0)_l$ of $B_0$ as in
Lemma~\ref{strat} such that
$\sigma|_{B^l_0}$ is a submersion onto a real-analytic submanifold
$A^l_0$ together with corresponding retractions $r^l\colon\Omega^l_0\to
A^l_0$ as in Step 1. Then each set
$$\2V^l:=\{(p,\L,t') \in \2M_{j_0+1}\times J^k_0\times (-1,1) :
\l_0(p,t')\in \Omega^l_0 , \,
(r^l(\l_0(p,t')),\L)\in B^l_0 \},$$
is a submanifold of $\2M_{j_0+1}\times J^k_0\times (-1,1)$.
As in the proof of Step 1 we see that
$\codim \2V^l\ge \codim W_0-1$ and hence
$\codim \2V^l > \dim_\R M+1$ in view of \eqref{codimw}.
We finally set
$$V^l:=\Pi^{-1}(\2V^l)\subset J^k(\2M_{j_0+1}\times(-1,1),\C^N),$$
where $\Pi\colon J^k(\2M_{j_0+1}\times(-1,1),\C^N) \to
\2M_{j_0+1}\times J^k_0\times (-1,1)$ is the natural projection.
Note that each $V^l$ is also a submanifold of
$J^k(\2M_{j_0+1}\times(-1,1),\C^N)$
of the same codimension as $V^l$.
By Thom's transversality theorem applied to $V^l$,
we obtain a residual set of maps $f\in C^{\infty}(\2M_{j_0+1}\times
(-1,1), \C^N)$ with $j^kf\cap V^l=\emptyset$ for all $l$, i.e.
$j^k_p\2\tau_{j_0}(\cdot,t')+t j^k_p f(\cdot,t')\notin W_p$
for all $t'\in (-1,1)$, $t\in (-1,1)\setminus\{0\}$, and $p\in \2M_{j_0+1}$.

We next choose a smooth function $\phi\in C^{\infty}(\2M_{j_0+1},\R)$
with compact support with $\phi(p)=1$ for $p$ in a neighborhood of
$\1M_{j_0+1}$. Furthermore, as a consequence of Property (2) for $j_0$,
the compactness of $\cup_{j\le j_0}
\1M_j$ and the fact that $W_p$ is closed in the open set in  $J^k_p$ of all
$k$-jets of germs of generic embeddings, we have the positive continuous distance
function
$$d(t):=\inf_{p\in \cup_{j\le j_0}
\1M_j}\dist(j^k_p\2\tau_{j_0}(\cdot,t), W_p) > 0, \, t\ne0,$$
where the distance is taken with respect to any fixed metric on
$J^k(M,X)$.
We now define
$$\2\tau_{j_0+1}(p,t):=
\begin{cases}
\2\tau_{j_0}(p,t) + \delta(t)\phi(p)f(p,t), & p\in \2M_{j_0+1},\\
\2\tau_{j_0}(p,t), & p\in M\setminus\2M_{j_0+1},
\end{cases}
$$
where $\delta(t)$ is a smooth real  valued function on $(-1,1)$ with
$\delta(0)=0$ and $\delta(t)>0$ for $t\ne 0$,
chosen sufficiently small such that
$$\dist\big(j^k_p\2\tau_{j_0+1}(\cdot,t) , j^k_p\2\tau_{j_0}(\cdot,t)
\big) < d(t), \quad t\in (-1,1)\setminus \{0\}, \quad p\in M,$$
and such that Properties (3) and (4) hold for $j=j_0+1$.
Here we use the standard fact that,  
for a given continuous positive function on $(-1,1) \setminus \{0\}$,
there exists a smooth function on $(-1,1)$, which is smaller than the given
function and positive away from $0$.
Then it follows that $\2\tau_{j_0+1}(p,t)$ satisfies Property (2) for
$j=j_0+1$. The remaining Properties (1) and (5) follow directly from the
construction.

We have now constructed the sequence $(\2\tau_j)_j$
satisfying Properties (1)--(5) above. Then we define the desired
deformation $\2\tau(p,t)$ as follows. Since the covering $(\2M_j)_j$
of $M$ is locally finite, there exists an integer function $\nu(p)\in\N$,
$p\in M$, such that every $p\in M$ has a neighborhood that is disjoint
from the union $\bigcup_{l\ge \nu(p)} \2M_l$. We set
$$\2\tau(p,t):=\2\tau_{\nu(p)}(p,t).$$
Then it is an immediate consequence of Property (5) that $\2\tau(p,t)$
is a smooth map from $M\times (-1,1)$ to $X$
and Properties (1), (2) and (4) imply the desired conclusions of
Theorems~\ref{strong} and ~\ref{strong-precise}.
\end{proof}

\section{Remarks and examples}\Label{remarks}

It should be mentioned that even a hypersurface in
$\C^2$, cannot, in general, be
approximated by, or deformed into, a 1-nondegenerate (i.e.\
Levi-nondegenerate) hypersurface.  Indeed, this is not possible if the
(scalar-valued) Levi form of the hypersurface changes sign. For
instance, the hypersurface $M\subset
\C^2$ given by $\{(z,w)\in \C^2: \Im w = (\Re z)^3\}$ cannot
be approximated by a Levi nondegenerate hypersurface, since
a Levi nondegenerate hypersurface would have either a strictly positive
or a strictly negative scalar Levi form (while the Levi form of $M$
takes both signs, depending on the sign of $\Re z$). However, it
follows from Theorems
\ref{weak-precise} and \ref{strong-precise} that any real hypersurface
in
$\C^2$ can be approximated by a 3-nondegenerate hypersurface. The
following example shows that in
$\C^2$  approximation by a 2-nondegenerate hypersurface is
impossible in general, so that these theorems are sharp for $N=2$ and
$m=3$.
\be\Label{no2} Consider the real-analytic hypersurface
$M\subset\C^2$ given by
\begin{equation}\Label{ge}
\Im w = 2\Re (z^3\bar z) + z\bar z\Re w.
\end{equation} It is easy to check that $M$ is
$3$-nondegenerate at $0$.  Let $W'\subset
J^3(\R^3,\C^2)$ be the set
of all
$3$-jets of
$2$-degenerate embeddings of $\R^3$ into $\C^2$.  As in
\eqref{harmfull}, we write the decomposition
  $$J^3_{0,0}(\R^3,\R) =(J^3)_h \oplus
(J^3)_{nh}.$$
We observe that $(J^3)_{nh}$ is a $4$-dimensional vector space and
that $W'$ is a submanifold of $J^3(\R^3,\C^2)$ of real
codimension
$3$.  Indeed, using Lemma \ref{split} with $k = 3$, $N=2$, and
$m=3$ (as well as the definition and invariance of $2$-degeneracy), we
see that, in the notation of the lemma, $W'$ is given by the vanishing of
the
$3$ real components of
$\Psi_4$ corresponding to $\psi_{z\1z}(0)$ and the real and imaginary
parts of $\psi_{z^2\1z}(0)$. We write $\Psi_4'$ for the mapping
corresponding to these three real components of $\Psi_4$. We now claim
that if
$\tau$ is the embedding from $\R^3$ to $\C^2$ whose image is $M$
(given by \eqref{ge}), the
$3$-dimensional submanifold $j^3\tau\subset J^3(\R^3,\C^2)$ intersects $W'$
transversally.  The fact that $M$ cannot be approximated by a
$2$-nondegenerate hypersurface will then follow from this
transversality.  Indeed, if  $\tau'$ is an embedding close to $\tau$, then
$j^3\tau'$ must also intersect $W'$ at some point, and hence
$\tau'(\R^3)$  cannot be  everywhere $2$-nondegenerate.

To prove the claimed transversality, we define the
mapping
$$
\R^3 \ni (x,y,s) \mapsto\mu(x,y, s):=\Psi_4'\big(j^3_{(x,y,s)}\tau\big).
$$
Since $W'$ is defined by the equation $\Psi_4'=0$, the transversality of
$j^3\tau$ to $W'$ is equivalent to the fact that $\mu$ is of rank $3$ at
the origin of $\R^3$.  The latter fact can be proved by a direct
calculation using the definition of $\Psi_4$ given in Lemma \ref{split},
and is left to the reader.
\end{Exa}

\br\Label{rem} Though a hypersurface $M\subset\C^2$
cannot be approximated in general by $2$-nondegenerate
ones, we claim  that such a
hypersurface can be approximated by a
$3$-nondegenerate one that is $2$-nondegenerate outside a
discrete subset. To see this, let $M$ be a $3$-dimensional manifold,
$W_3'\subset J^4(M,\C^2)$  the set
of all
$4$-jets of
$3$-degenerate embeddings of $M$ into $\C^2$, and $W_2'\subset
J^3(M,\C^2)$  the set
of all
$3$-jets of
$2$-degenerate embeddings of $M$ into $\C^2$. Note that $\codim
W_3' = 6$ and $\codim W_2' = 3$.  By Thom transversality theorem, the
set
\begin{equation}\Label{fork} \{ f \in C^\infty(M,\C^2): j^3f \pitchfork
W_2', \ j^4f\cap W_3'=\emptyset\}
\end{equation}
is residual, and hence any embedding
   $\tau:M \to \C^2$ can be approximated by an embedding in
\eqref{fork}.  The rest of the claim follows from the fact that, by
dimension, if $j^3f $ is transversal to
$W_2'$, then the intersection of the two sets is necessarily discrete (or
empty).
\er

In particular, it follows from Remark~\ref{rem}  that any compact
hypersurface in $\C^2$ can be approximated by another one which is
$3$-nondegenerate everywhere and $2$-nondegenerate except at finitely
many points.  Our final example shows that for any positive integer
$\ell$, there is a real-analytic boundary of a
bounded domain in $\C^2$ that cannot be approximated by a
hypersurface that is $2$-nondegenerate except at $\ell-1$ points.

\begin{Exa} Let $\ell$ be a positive integer, and let
$(r,\omega)\in\R_+\times S^3$ be the usual polar coordinates in
$\C^2 \cong \R^4$.  We shall show that there a real polynomial $P$
defined in $\R^4$ such that the boundary of the domain
\begin{equation}\Label{D}
D:=\{(r,\omega) \in \C^2: 0\le r<e^{P(\omega)}\}
\end{equation}
satisfies the above property.  For this, we choose any $\ell$ distinct
points
$\{p_1,\ldots, p_\ell\}$ on the unit sphere $S^3\subset \C^2$ and let
$M$ be the hypersurface given in Example \ref{no2}. The reader can
easily check that there is a real polynomial $P$ defined in $\C^2$
vanishing at $p_1,\ldots,p_{\ell}$ such that the following hold:
\begin{enumerate}
\item If $D$ is given by \eqref{D} then its boundary, $\partial D$, is
tangent to the  (real) tangent hyperplane to
the sphere at $p_k$, for $k=1,\ldots, \ell$.

\item For each $k$,
  $k=1,\ldots,\ell$, the $4$-th jet at $p_k$ of the hypersurface $\partial D$
is the same as that of
$M$ at $0$, after an appropriate translation and complex linear
transformation in $\C^2$.

\end{enumerate}
Note that the conditions above involve only the derivatives of $P$ at
$p_k$, $k=1,\ldots,\ell$, of order $\le 4$, and that $P$ is obviously not unique.

We now observe that, by the transversality argument used
in Remark~\ref{rem}, any real hypersurface
$\2M\subset\C^2$ that is tangent to $M$ (given in Example~\ref{ge})
at the origin of order at least $4$, has the same property as $M$, i.e.\
$\2M$ cannot be approximated by smooth $2$-nondegenerate
hypersurfaces. We conclude  that any smooth compact hypersurface
sufficiently close to $\partial D$ must have at least $\ell$
$2$-degenerate points.   Recall  that the approximating hypersurface
can always be chosen to be everywhere $3$-nondegenerate as a
consequence of Theorem~\ref{weak-precise}.

\end{Exa}

\end{document}